\documentclass[12pt,eqno]{article}
\usepackage{latexsym}
\usepackage{amsfonts,amsmath,amssymb}                                                                                                                                                                                                                                                                                                                                                                                                                                                                                                                                                                                      
\usepackage{amsthm}
\newtheorem{theorem}{Theorem}[section]
\newtheorem{lemma}{Lemma}[section]

\newtheorem{corollary}{Corollary}[section]

\usepackage[pdftex]{graphicx}
\usepackage[latin1]{inputenc}
\usepackage{epstopdf}
\usepackage{setspace}
\usepackage{tikz}
\usepackage[pagewise]{lineno}
\usepackage[toc]{appendix}

\numberwithin{equation}{section}
\numberwithin{theorem}{section}

\newtheorem*{prth1.1}{Proof of Theorem 1.1}
\newtheorem*{prth1.2}{Proof of Theorem 1.2}
\newtheorem*{proofcor}{Proof of Corollary 1.1}
\numberwithin{proposition}{section}
\numberwithin{lemma}{section}
\numberwithin{remark}{section}
\numberwithin{equation}{section}
\begin{document}
\doublespacing
\setcounter{secnumdepth}{3}
%\numberwithin{equation}{section}
\newcommand{\bbox}{\vrule height.6em width.6em
depth0em} %%%%% Black Box
\parindent0in
\def\lhead{Behaviour in time of solutions}
\thispagestyle{empty}
\setcounter{page}{1}
\noindent
\thispagestyle{empty}
\setcounter{page}{1}
\noindent\\
\thispagestyle{empty}
\setcounter{page}{1}
\noindent
\begin{center}
{ \bf{
Behaviour in time of solutions to  fourth-order parabolic systems with time dependent coefficients.}}\\
\vspace*{1cm}
M.Marras \footnote{ Dipartimento di Matematica e Informatica, Universit\'a di Cagliari, via Ospedale 72, 09124 Cagliari (Italy), mmarras@unica.it
 }
 and S. Vernier Piro \footnote
{Facolt\'a di Ingegneria e Architettura, Universit\'a di Cagliari,
 Via Marengo 2, 09123 Cagliari (Italy), svernier@unica.it}  \\
\end{center}
\vspace*{0.5cm}

%%==========================================================%%
%%==============                                                    ==============%%
%%======                                                                                    ======%%
%%====                                                                                        ====%%
%%==                               ABSTRACT                                      ==%%
%%===                                                                                           ====%%
%%======                                                                                 ======%%
%%==============                                                     ==============%%
%%==========================================================%%

{\bf Abstract.}
This paper deals with a class of initial-boundary value problems for  nonlinear  fourth order  parabolic systems  with time dependent coefficients  in a bounded domain $\Omega \subset \mathbb{R}^N, N\geq 2$. 
Introducing suitable conditions on the source terms, we obtain a time interval   $[0,T],$ where the solution remains  bounded by  deriving  a lower bound  $T$ of $t^*$. Moreover,
we establish conditions on the shape of the spatial domain  and on data sufficient to guarantee that the solution blows up in finite time $t^*$, deriving an upper bound for $t^*$.
%%%%%%%%%%%%%%

\vskip1.truecm
%%==========================================================%%
%%==============                                                    ==============%%
%%======                                 Section 1                                   ======%%
%%====                                                                                        ====%%
%%==                                                                                              ==%%
%%===                                      Introduction                                 ====%%
%%======                                                                                 ======%%
%%==============                                                     ==============%%
%%==========================================================%%

\section {Introduction}
We deal with the following  fourth-order parabolic system with time dependent coefficients
\begin{equation}
\left\{ \begin{array}{l}
\label{sys}
u_t  +  \delta_1(t)\Delta^2 u-  h_1(t) \Delta u = k_1(t)  f_1(v),  \quad
 x  \in \Omega , \ \  t \in (0,t^*),\\
v_t  +  \delta_2(t) \Delta^2 v - h_2(t) \Delta v= k_2(t)  f_2(u),  \quad
 x  \in \Omega , \ \  t \in (0,t^*),\\
 \end{array} \right.
\end{equation}
\begin{equation}
\label{dir} 
 (u,v)=0,   \quad  \left( \frac {\partial  u}{\partial n} ,  \frac {\partial  v}{\partial n} \right)=0,    \quad {x}  \in \partial \Omega , \ \  t \in (0,t^*),
\end{equation}
\begin{equation}
\label{initial}
(u ({x},0), v ({ x},0)) = (u_0 ({x}), v_0({x})),  \quad { x} \in \Omega,    \\
\end{equation}
where $\Omega$ is a bounded domain in $\mathbb{R}^N, \, N\geq 2$,  with smooth boundary $ \partial \Omega$,    $\displaystyle \frac {\partial  }{\partial n}$  is the outward normal derivative on the boundary $\partial \Omega$, $f_1$ and $f_2$ are non negative functions, \ 
  $t^*$ is the  blow-up time of the solution; $(u_0(x),v_0(x))$ are non negative functions in $\Omega$, $\delta_i, k_i, h_i \ i=1,2,$ are   in general positive bounded functions of $t$.\\
In the literature the great part of the results concerns the case of only one equation and this because  fourth order parabolic equations describe  a variety of physical processes (see \cite{GMP}). In particular, in the thin film theory (see \cite  {H}, \cite {KSW}), \cite{OR}) such problems describe the evolution of the epitaxial growth of thin film and the following equation is introduced:
\begin{equation}
\label{thin-film}
u_t  + \Delta^2 u - A_1 \Delta u- A_2 div( |\nabla u|^2 \nabla u)= g(x,t,u)
\end{equation}
where $ u(x,t)$ is  the height from the surface of the film, $\Delta^2 u$ denotes the capillarity-driven surface diffusion, $A_1\Delta u$ describes  the diffusion due to 
\\
evaporation-condensation, $A_2 \ div( |\nabla u|^2 \nabla u)$ represents the upward hopping of atoms, $g$ corresponds  to the mean  deposition flux.
In  \cite {KSW}  King, Stein and 
Winkler  showed  for the solutions of \eqref{thin-film} existence,  uniqueness  and  regularity in an appropriate function space.
   Xu et al. in \cite{XU} investigated  the equation \eqref{thin-film} with $A_2=0$ and $g(x,t,u)=g(u)$ in a bounded domain in ${\mathbb R}^N$ and, by using the potential well method,  
showed that the solutions exist globally or blow-up in finite time, depending on whether or not the initial data  are in the potential well.
When the Lapacian term in \eqref{thin-film}  is replaced by p-Laplacian and $g(u)$ is of power type, with a modified potential well method, Han in \cite{H}  obtained  again global existence and finite time blow-up  for the solutions  when the initial energy is subcritical and critical.
By  different methods Philippin in \cite{P},  with  $A_1=A_2=0, \ g(u)=k(t) |u|^{p-1} u, $  proved that the solution  cannot exist for all time, i.e. it  blows up in $L^2$-norm and an upper
bound for $t^*$ is derived; 
in addition, he  constructed (under certain conditions on the data) a lower bound
for $t^*$ by using  a first order differential inequality method.
 Escudero,  Gazzola and  Peral in  \cite{EGP}  proved existence and blow-up  results for the solutions of the equation
$$u_t + \Delta^2 u = det( D^2 u) + \lambda h(x,t), $$
under  Dirichlet  boundary  conditions,
which models epitaxial growth processes and where  the evolution is dictated by the competition between  the determinant of the Hessian matrix of the solution and the bilaplacian (see also  \cite{EGHPT}). 
\\
%%%%%   Winkler 2011
Recently Winkler in  \cite{W} considered  the equation
$$
u_t  = - \Delta^2 u - \mu  \Delta u- \lambda \Delta( |\nabla u|^2) + f(x)
$$
 when $\Omega$ is a bounded convex domain in $\mathbb{R}^N$, under the conditions $\displaystyle  \frac {\partial  u}{\partial n} = \frac {\partial  \Delta u}{\partial n} =  0$
on the boundary $\partial \Omega$ with bounded initial data, and proved the existence of  global weak  solutions with spatial dimension $N\leq 3$, under suitable conditions on data.\\

Less attention was given to the case of parabolic fourth order systems. In (\cite{PVp}, Sec.4) Philippin and Vernier Piro 
investigated  a simplified version of the system \eqref{sys}- \eqref{initial}  i.e. with $\delta_i (t)=1, h_i(t)=0,$
and   $\Omega$  a bounded domain in $\mathbb{R}^N$, with $N=2$ or $N=3$.  They  proved that there exists a safe interval of existence $[0,T]$, $T$ a lower bound of $t^*,$ extending a method used in the study of second order parabolic systems (see i.e. \cite{MV1}, \cite{MV2}, \cite{PP}). The method is based on introducing suitable functionals which satisfy   first order differential inequalities  from which upper and/or  lower bounds of the blow-up time are derived.\\

The  aim of this article is  to investigate the behavior  near the blow-up time of the solutions in the case of a more general  system as (\ref{sys})--(\ref{initial}),  where  the dimension  of the spatial domain $\Omega$ is $N \geq 2$: we  introduce conditions on data sufficient to prove that the solution must blow up in finite time $t^*$, deriving an upper bound of $t^*$;
moreover,  we prove that there exists an interval where the solution remains bounded, by deriving a lower bound of $t^*$. 
 
In this contest, we consider  non negative solutions of the system \eqref{sys} - \eqref{initial}, motivated by our aim to  investigate the behavior of the solutions which approach to $+\infty$ as $t$ approaches the  finite blow-up time  $t^*$. For the definition of solution to  the system \eqref{sys}- \eqref{initial} we extend the definition  in \cite{EGP}, Theorem 3.1: for  some $T>0$ and $f_1$ and $f_2 \in L^2 (0,T; L^2 (\Omega))$, $u_0$ and $v_0 \in W_0^{2,2}(\Omega)$,  the functions $u$ and $v$ belong to the space  $ C(0,T; W_0^{2,2}(\Omega)) \cap L^2(0,T; W^{4,2}(\Omega)) \cap 
W^{1,2}  (0,T; L^2 (\Omega))$.

Throughout the paper we assume $f_1$ and $f_2 \in L^2 (0,T; L^2 (\Omega))$, $u_0$ and $v_0 \in W_0^{2,2}(\Omega)$.
 We denote by $||\cdot||_ {\sf r} $ the $L^{\sf r} (\Omega)$ norm for $1\leq {\sf r} \leq \infty$ and  by $||u||_{W^{2,2}_0(\Omega)}= ||\Delta u||_2=\left( \int_{\Omega} (\Delta u)^2 dx \right )^{\frac 12 }$  the norm in $W^{2,2}_0(\Omega) $.
 
 The first investigation in this paper is to determinate a safe time interval of existence of the solution $(u,v)$ of \eqref{sys}-\eqref{initial}, say $[0,T]$, with $T$ a lower bound of $t^*$. We remark that in this investigation $\Omega$ is not assumed to be a ball.
With the aim to obtain an interval where the solution remains bounded, we introduce the functional
\begin{equation} \label{Phi}
\Phi(t)=||\Delta u||^2_2+||\Delta v||^2_2=\Phi_1(t)+\Phi_2(t) ,
\end{equation}
with initial value 
\begin{equation} \label{Phi_0}
 \Phi(0):=\Phi_0= ||\Delta u_0||^2_2+||\Delta v_0||^2_2.
 \end{equation}

%%% lower bound
\begin{theorem}  (Lower Bound).
\label{lower}
Let $\Omega $ be a bounded domain in ${\mathbb R}^N$. Let $(u,v)$ be a non negative solution of the system \eqref{sys}-\eqref{initial}  and let $\Phi(t)$ be defined in \eqref{Phi}-\eqref{Phi_0}. Assume 
 
\begin{equation}
\label{f_1 f_2 lower}
f_1(s) \leq s^p, \ \ \ f_2(s) \leq s^q, \ \ \ s>0, \   \ with \ p>q>1,
\end{equation}
$p, q$ satisfying the condition \eqref{Sob} in Lemma \ref{lemma1} below.
 Then  $\Phi(t)$ remains bounded in the interval $[0,T],$ with $T<t^*$  and
 \begin{equation}
 \label{Ti}
\begin{cases}
T:= \frac{1}{B(p-1)} \ln\Big[1 +\frac {B}{2A} \Phi^{1-p}_0\Big] ,  & if \ \  p>q, \ \ A, \, B >0, \\[2mm]
T:= \frac{1}{B(p-1)} \ln\Big[1 +\frac {B}{\widetilde A} \Phi^{1-p}_0\Big], & if \ \ q=p \ \ \widetilde A, \, B >0.\\[2mm]
\end{cases}
\end{equation}
 \end{theorem}

To obtain an upper bound of the blow-up time
we introduce the  functional
\begin{equation}
\label{Psi}
\Psi(t)= \int_{\Omega} u\phi_1 dx+ \int_{\Omega} v\phi_1 dx = \Psi_1(t) + \Psi_2(t),  
\end{equation}
with
\begin{equation}
\label{Psi0}
 \Psi(0)=\Psi_0= \int_{\Omega} u_0\phi_1 dx+ \int_{\Omega} v_0 \phi_1 dx =\Psi_1(0) + \Psi_2(0)
\end{equation}

and $\phi_1$ the first eigenfunction associate to the first eigenvalue $\Lambda_1$ of the biharmonic eigenvalue problem with Dirichlet boundary conditions:

\begin{equation}
\label{eigen}
 \Delta^2 \phi   = \Lambda \phi , \quad x\in \Omega \subset \mathbb{R}^N, \ N\geq 2, 
 \end{equation}
  \begin{equation}
\label{direll} 
 \phi=0, \   \quad   \frac {\partial  \phi}{\partial n} =0 ,   \quad {x}  \in \partial \Omega , 
  \end{equation}
with $\phi$  normalized by
\begin{equation}
\label{normalize} 
||\phi|| _2^2=1.
\end{equation}
 For all $\phi \neq 0$, $\Lambda_1$ satisfies 
 \begin{equation}
\label{eigen2}
||\phi||_2^2 \leq    \Lambda_1^{-1} ||\Delta \phi||^2_2,
\end{equation}
(see \cite{AH}).
The problem \eqref{eigen} is closed related to the biharmonic differential
equation 
$$
 \Delta^2 \phi   = f
$$
with the same boundary conditions \eqref{direll}, which describes characteristic  vibrations of a clamped plate.
 It is well-known that the biharmonic operator under Dirichlet boundary conditions does not satisfy the positivity preserving property in general domains; nevertheless  it holds when  $\Omega=B_R$ the ball in $\mathbb{R}^N$ with radius $R$. 
 We remark that the balls are not the only sets where the property holds.
 (See e.g. \cite{GGS},  \cite{GS} and \cite{S}).
As a consequence the functional $\Psi(t)$  defined in \eqref{Psi} and \eqref{Psi0} is non negative.\\

 We prove the following 

\begin{theorem}(Blow up and upper bound).\\
\label{BlowUp}
Let $\Omega = B_R$ be the $N$-dimensional ball, $N \geq 2$. Let $(u,v)$ be a non negative solution of  the system \eqref{sys}-
\eqref{initial} with $h_i=0$
and let  $\Psi(t)$ be defined in \eqref{Psi}-\eqref{Psi0}. Assume 

\begin{equation}\label{f_1 f_2}
f_1(s) \geq s^p, \ \ \ f_2(s) \geq s^q, \ \ \ s>0, \   \ with \ p>q>1,
\end{equation}
and the initial data are so large that
\begin{align} \label{H(Psi)}
H(\Psi_0)>0, \ \  \ H(\Psi(t)):=-\Lambda_1\delta \Psi + 2^{1-q} c \, \Psi^q  - c\, Q, \ \ c,Q>0
 \end{align} 
with $\Lambda_1$ the first eigenvalue of the  problem \eqref{eigen}- \eqref{normalize}, then $\Psi$  blows up in a finite time $t^*$  in the sense that
\begin{equation*}
\lim_{t\rightarrow t^*} \Psi(t)= \infty,
\end{equation*}
with the upper bound  
\begin{equation}\label{upper}
T_0=  \int_{\Psi_0}^{\infty}  \frac {d\eta}{H(\eta)}\geq t^*.
 \end{equation}
\end{theorem}
In the particular case $q=p,$ we achieve the following corollary.
\begin{corollary} \label{corollary}
Under the hypotheses of Theorem \ref{BlowUp}, if $q=p$ and if there exist two positive costants $\delta, \bar c$ such that the initial data $(u_0,v_0)$  satisfy the following condition
  \begin{equation} \label{Psi0 p=q}
\Psi_0> \left( \frac { \delta \Lambda_1}{\bar c}\right)^ {\frac 1{p-1}  }, 
 \end{equation}

then $\Psi(t)$ must blow up at time $t^* \leq \overline T $ with 
\begin{equation} \label{bar T, p=q}
\overline T:=- \frac 1 {(p-1)\delta \Lambda_1} \log \Big\{  1-  \frac {\delta  \Lambda_1} {\bar c \Psi_0^{p-1}}  \Big\}.
\end{equation}
\end{corollary}

%%%%%%%%%%%
%%%%%%%%%%%
%%%SECTIONS

The scheme of this paper is the following:
in Section \ref{lowerB} we prove Theorem \ref{lower} showing   that the functional  $  \Phi(t) $ remains bounded on some time interval $(0,T),$ where in some  particular cases  $T$ may be explicitly computed in terms of the data of the system (\ref{sys}). Clearly this value of $T$ provides a lower bound for blow-up time $t^*$ of the solution $(u,v).$
As a consequence,  we have that also the  $L^2$-norm of the solution remains bounded in $(0,T),$   since from the eigenvalues theory (see \cite{H}) we have  
$$||u||_2^2 + ||v||_2^2  \leq  \Lambda_1^{-1}( ||\Delta u||^2_2+||\Delta v||^2_2), $$
 where $\Lambda_1$ is the first eigenvalue of the  problem \eqref{eigen}-\eqref{normalize} applied to $u$ and $v$.\\
In Section \ref{upperB}  we consider the system (\ref{sys})- (\ref {initial}) with $h_i(t)=0$, and  the spatial domain  the N-dimensional ball $B_R(0)$  and we prove Theorem \ref{BlowUp} and Corollary \eqref{corollary}. \\
%==========================================================%%
%%==============                                                    ==============%%
%%======                                 Section 2                                   ======%%
%%====                                                                                        ====%%
%%==                                                                                              ==%%
%%===                                     Preliminaries.                               ====%%
%%======                                                                                 ======%%
%%==============                                                     ==============%%
%%==========================================================%%

%%%%%%%%%%%%%%
  
%%%%%%%%%%%%%

\section{ Lower bound of $t^*$ }\label{lowerB}
In this section we consider the system \eqref{sys}-\eqref{initial}  and we assume that the solution blows up at finite time $t^*$ in the sense that $\displaystyle \lim_{t \rightarrow t^*} \Phi(t)= \infty$ with $\Phi$ defined in \eqref{Phi}-\eqref{Phi_0}.
 We state now  a known result which will be needed in this section.\\
From the Rellich-Kondrachov theorem on the embedding $W^{{\mathsf p},m} \subset L^{\sf r}$ (see \cite{GGS} Theorem 2.4), the following Lemma is derived for $\sf p=2.$
%\subsection{A Sobolev-Poincar$\acute{e}$ inequality}
 \begin{lemma}\label{lemma1}
 Let $\Omega$ be a bounded domain in $  \mathbb{R}^N$. Let $\sf r$ be an arbitrary number  with
 $ 2\leq {\sf r} <+ \infty,\ if \  N< 4$ and  $ 2\leq {\sf r} < \frac{2N } {N-4 }, \ if  \ N>4.$
 Then for any $w\in W_0^{2,2}(\Omega)$ there exists a constant $S=S({\sf r},\Omega)$ such that
\begin{equation}
\label{Sob}
|| w||_{\sf r} \leq S ||\Delta w||_2.
\end{equation} 
\end{lemma}
 To derive a lower bound of the blow up time $t^*$ of the solution of \eqref{sys}-\eqref{initial}, we construct a first order differential inequality for $\Phi$  defined on \eqref{Phi} and prove the Theorem \ref{lower}.
\begin{prth1.2}
 {\rm 
 Testing the first equation in \eqref{sys} with $ \Delta\Delta u $ we have
\begin{equation}
\label{test}
\int_{\Omega} u_t  \Delta^2 u+\delta_1 \int_{\Omega}(\Delta^2 u)^2 - h_1
\int_{\Omega}\Delta u \  \Delta^2u =  k_1\int_{\Omega} f_1\Delta^2 u.
\end{equation}

By using two times the $\epsilon-$Young inequality we get
\begin{equation}
\label{test1}
 \int_{\Omega} f_1\Delta^2 u
\leq  \frac{  \epsilon_1}2 \int_{\Omega} f_1^2 + \frac 1 {2 \epsilon_1 }  \int_{\Omega} (\Delta^2 u)^2,
\end{equation}
and
\begin{equation}
\label{test2}
\int_{\Omega}\Delta u \  \Delta^2 u\leq   \frac{  \epsilon_2}2 \int_{\Omega} |\Delta u|^2 + \frac 1 {2 \epsilon_2 }  \int_{\Omega} (\Delta^2 u)^2
\end{equation}
with $ \epsilon_1=\epsilon_1(t)$  and $ \epsilon_2=\epsilon_2(t)$
  two arbitrary  positive functions.\\

Taking into account \eqref{f_1 f_2 lower} we get
\begin{align}\label{I3}
\frac{  \epsilon_1}2 \int_{\Omega} f_1^2 + \frac 1 {2 \epsilon_1 }  \int_{\Omega} (\Delta^2 u)^2dx \leq  \frac{ \epsilon_1 }2
  || v^{p}||_2^2  
 + \frac 1{ 2 \epsilon_1}    ||\Delta^2 u||_2^2  \Big),
\end{align}

Inserting \eqref{test1} , \eqref{test2} and \eqref{I3}     in \eqref{test}, it follows

\begin{align}
\label{test3}
& \frac 12 \ \frac {d}{dt} \int | \Delta u|^2\\
&  \leq -
 \delta_1 || \Delta^2 u||_2^2  +
 \frac{ h_1 \epsilon_2}2 ||\Delta u||_2^2 + \frac {h_1} {2 \epsilon_2 }  ||\Delta^2 u||_2^2+ \frac{k_1  \epsilon_1}2  || v^{p}||_2^2   + \frac{k_1} {2 \epsilon_1 } ||\Delta^2 u||_2^2.
\end{align}
 
Now  using \eqref{Sob} in Lemma \ref{lemma1} with ${\sf r} = 2 p$, we obtain
\begin{equation}
\label{v^2p}
|| v^{p}||_2^2 \leq  S^{2p} \Big(   ||\Delta v||_2^{2} \Big)^p.
\end{equation}
At last  we can write
\begin{align}
\label{Phi'1}
 \Phi_1'(t)\leq &  \Big(\frac{ h_1}{ \epsilon_2} + \frac{k_1}{ \epsilon_1} - 2 \delta_1 \Big)
  ||\Delta^2 u||^2  +  k_1\epsilon_1 S^{2p} \Big( ||\Delta v||_2^{2} \Big)^p + h_1 \epsilon_2  ||\Delta u||_2^2 
\end{align}
In the same way we derive
\begin{align}\label{Phi'2}
 \Phi_2'(t) \leq &  (\frac{h_2}{ \epsilon_4} +\frac{ k_2}{ \epsilon_3} - 2 \delta_2)   ||\Delta^2 v||^2_2 
+k_2  \epsilon_3 S^{2q}\Big(  || \Delta u||_2^2 \Big)^q      
  +  h_2 \epsilon_4   ||\Delta v||_2^2
\end{align}
with $ \epsilon_3=\epsilon_1(t)$  and $ \epsilon_4=\epsilon_2(t)$
  two arbitrary  positive functions.
Since $\Phi'=\Phi'_1+ \Phi_2'$, using \eqref{Phi'1} and \eqref{Phi'2} and choosing the functions $\epsilon_i, \ i=1,..,4$ such that the first terms in \eqref{Phi'1} and \eqref{Phi'2} vanish, we arrive at
\begin{align}
\label{Phi' bis}
 \Phi'(t)\leq &  k_1  \epsilon_2 S^{2p} \Big(  ||\Delta v||_2^2\Big)^p \!\!+ k_2  \epsilon_3 S^{2q}\Big(  ||\Delta u||_2^2 \Big)^q \!\!
+ h_1 \epsilon_2  ||\Delta u||_2^2 + h_2 \epsilon_4   ||\Delta v||_2^2 \\
\notag = & k_2  \epsilon_3 S^{2q} \Phi_1^q +  k_1  \epsilon_2 S^{2p}  \Phi_2^p +h_1 \epsilon_2 \Phi_1 + h_2\epsilon_4 \Phi_2\\
 \notag \leq &    k_2  \epsilon_3 S^{2q} \Phi^q +  k_1  \epsilon_2 S^{2p}  \Phi^p + B  \Phi,
\end{align}
with $ \displaystyle B:= \max \{B_1, B_2 \}$ with
 $B_1=\sup_t \{h_1(t) \epsilon_2(t) \}$ and $B_2= \sup_t\{ h_2(t) \epsilon_4(t) \}.$
Since $\Phi(t)$ blows up in a finite time $t^*,$ then there exists a time $t_1 >0$ such that $\Phi(t) \geq \Phi(t_1)=\Phi_0$ for all $t \geq t_1$ and we can write with $p>q$
\begin{equation} \label{Phi(t_0)}
\Phi^q (t) \leq \Phi^p(t)  \ \Phi^{q-p}_0, \ \ \ \forall \ t\in (t_1, t^*).
\end{equation}
In view of \eqref{Phi(t_0)}, the inequality \eqref{Phi' bis} becomes
\begin{align}
\label{Phi' final}
 \Phi'(t)\leq   2A \Phi^p + B \Phi,
\end{align}
with $\displaystyle A:= \max  \{A_1,A_2 \} $ with  $A_1= \sup_t  \{ k_2(t) \epsilon_3(t) S^{2q}  \Phi^{q-p}_0  \}$ and \\ $A_2= \sup_t  \{ k_1(t)  \epsilon_2(t) S^{2p}\}.$ 
Integrating \eqref{Phi' final} we arrive at:
\begin{align*}
e^{B(p-1) t} \Phi^{1-p} (t) - \Phi_0^{1-p} &  \geq - 2A(p-1) \int_{0}^{t}  e^{B(p-1) \tau} d\tau	\\
 &= 
 \frac{-2A}{B(p-1)} \Big[ e^{B(p-1) t} -1\Big].
\end{align*}

Finally,  letting  $t\rightarrow t^*$, we obtain a lower bound $T$ of the blow-up time with:
\begin{equation*}\label{T}
 T:= \frac{1}{B(p-1)} \ln\Big[1 +\frac {B}{2A} \Phi^{1-p}_0\Big]\leq t^*.
\end{equation*}

In the particular case $p=q$, we replace \eqref{Phi' bis}  with
\begin{equation*}\label{p=q}
\Phi' \leq   k_2 \epsilon_3 S^{2p}  \Phi_1^p +  k_1  \epsilon_2  \Phi_2^p + B \Phi
\end{equation*}
and  by using the basic inequality $a^{\gamma} + b^{\gamma} \leq (a+b)^{\gamma},$ $\gamma >1, \ \ a,b >0$, we get
\begin{equation}\label{Phi' final p=q}
 \Phi' \leq   \widetilde A \Phi^p + B \Phi,
 \end{equation}
   with $\displaystyle \widetilde A =S^{2p} 
   \max  \{\widetilde A_1, \widetilde A_2\} $ with  $\widetilde A_1= \{k_2(t) \epsilon_3(t)   \Phi^{q-p}_0\}$ and \\ $\widetilde A_2= \sup_t \{k_1(t)  \epsilon_2(t)\}.$ 
\\
Integrating \eqref{Phi' final p=q}, we lead to the lower bound
$$
T:= \frac{1}{B(p-1)} \ln\Big[1 +\frac {B}{\widetilde A} \Phi^{1-p}_0\Big] \leq  t^*.
$$
We have obtained  \eqref{Ti} and the Theorem \ref{lower} is proved.
\qed}
\end{prth1.2}

{\bf Remark 1}\\
The existence of a lower bound $T$  for  the blow-up time to the functional $\Phi(t)$ has a consequence  that the interval $[0,T]$ is a safe interval of existence of the solution,  since, as  observed in the introduction,
$$||u||_2^2 + ||v||_2^2  \leq  \Lambda_1^{-1}( ||\Delta u||^2_2+||\Delta v||^2_2).$$
{\bf Remark 2}\\
We may obtain a further  lower bound of the blow up time $t^*$ (easier to be computed) by estimating $\Phi(t)$ on the right of  \eqref{Phi' final p=q}  in the following way:  arguing as in \eqref{Phi(t_0)} with $q=1$ we get
\begin{equation} \label{Phi(t_1)}
\Phi (t) \leq \Phi^p(t) \  \Phi^{1-p}_0, \ \ \ \forall \ t\in (t_1, t^*).
\end{equation}
By replacing  \eqref{Phi(t_1)} in \eqref{Phi' final} we arrive at
\begin{equation} \label{Phi' remark}
\Phi'(t) \leq  \ K\Phi^p,
\end{equation}
with $K:= 2A + B \Phi^{1-p }_0.$
Integrating \eqref{Phi' remark}  in the interval $(t_1, t)$ and letting $t\rightarrow t^*$, we have
\begin{equation*}
\int_{\Phi_0}^{\infty} \frac {d\eta}{\eta^p} = \int_{\Phi(t_1)}^{\infty}  \frac {d\eta}{\eta^p} \leq K  \int_{t_1}^{t^*} d \tau \leq K \int_0^{t^*}  d\tau = K t^*,
\end{equation*}
 from which we obtain the following lower bound
\begin{equation*}
\widetilde T= \frac{\Phi^{1-p}_0}{ (p-1)K} \leq t^*.
\end{equation*}

%%==========================================================%%
%%==============                                                    ==============%%
%%======                                 Section 3                                   ======%%
%%====                                                                                        ====%%
%%==                                                                                              ==%%
%%===                                    Upper bound.                               ====%%
%%======                                                                                 ======%%
%%==============                                                     ==============%%
%%==========================================================%%

\section{Upper bound of the blow up  time }\label{upperB}
In this section the system \eqref{sys}-\eqref{initial}  is simplified by the assumption $h_i=0$. We recall that, as already observed in the introduction,  the functional $\Psi(t)$ defined in \eqref{Psi}-\eqref{Psi0}, is non negative.

%%%%%%%%%%%%%%%%%%% 

We derive sufficient conditions on the data which guarantee that the solution $(u,v)$ blows up in finite time.\\ 
\begin{prth1.1}
{\rm 
Testing the first equation with $h_1 =0$ in \eqref{sys}  with $\phi_1$,  we have
\begin{equation}
\label{u_t}
\int_{\Omega} u_t \phi_1 dx+\delta_1  \int_{\Omega}(\Delta^2 u) \phi_1 dx = k_1\int_{\Omega} f_1 \phi_1 dx.
\end{equation}
For the first term on the left in \eqref{u_t}, with $\Psi_1$ defined in \eqref{Psi},
$$
\int_{\Omega} u_t \phi_1 = \frac{ d} {dt }\int_{\Omega} u \phi_1 = \Psi_1'(t)
$$
 since $\phi_1 $ doesn't depend on t. Then
 $$
 \frac{ d} {dt }\int_{\Omega} u \phi_1 dx= - \delta_1 \int_{\Omega}(\Delta^2 u) \phi_1 dx + k_1\int_{\Omega} f_1 \phi_1 dx.
 $$

Using  the hypothesis \eqref{f_1 f_2} we obtain 
\begin{equation}\label{Psi'_1}
\Psi'_1 = \int_{\Omega} (k_1 f_1- \delta_1\Delta^2 u)\phi_1dx \geq k_1 \int_{\Omega} v^p \phi_1 dx- \delta_1 \int_{\Omega}  \Delta^2 u \  \phi_1 dx.
 \end{equation} 
In the first term of \eqref{Psi'_1} we use the Holder's inequality and \eqref{normalize}  leading to\\
\begin{equation}\label{vp_phi1}
 \int_{\Omega} v^p\phi_1dx  \geq |\Omega|^{-\frac {p-1} 2 }\Psi_2^p.
 \end{equation}
 Using the second Green's identity and \eqref{eigen}, the second term of \eqref{Psi'_1} can be estimate as follows
\begin{equation}\label{DDelta_u_phi}
- \int_{\Omega} \phi_1 \Delta^2 u \  dx = -  \int_{\Omega} u  \Delta^2 \phi_1  dx = -\Lambda_1\int_{\Omega}  u \phi_1 dx = -\Lambda_1 \Psi_1.
 \end{equation}
Substituting \eqref{vp_phi1} and \eqref{DDelta_u_phi} in \eqref{Psi'_1} we arrive at
\begin{equation}\label{Psi'_1bis}
\Psi'_1(t) \geq  k_1 |\Omega|^{-\frac {p-1} 2 }\Psi_2^p (t)-\delta_1(t)\Lambda_1 \Psi_1 (t)=c_1(t) \Psi_2^{p} -\delta_1(t) \Lambda_1 \Psi_1, 
 \end{equation}
  with $c_1(t):= k_1(t) |\Omega|^{-\frac {p-1} 2 }$.\\
  
 Testing the second  equation with $h_2 =0$ in \eqref{sys}  with $\phi_2$,  
 a  computation similar to the  previous leads to
  \begin{equation}\label{Psi'_2}
\Psi'_2 (t)\geq c_2(t) \Psi_1^q - \delta_2(t) \Lambda_1 \Psi_2,
 \end{equation}
with $c_2(t):= k_2(t)|\Omega|^{-\frac {q-1} 2 }$.\\
Adding \eqref{Psi'_1bis} and \eqref{Psi'_2}, we have
\begin{equation}\label{Psi'bis}
\Psi'  \geq  - \delta\Lambda_1(\Psi_1 + \Psi_2) + c(\Psi_2^p + \Psi_1^q),
 \end{equation}
  where $\displaystyle c:= 
  \min\{C_1, C_2 \},$ with $C_1= \inf_t\{c_1(t)  \}$ and $C_2= \inf_t\{c_2(t)  \}$
   and $\displaystyle \delta:= \max \{D_1, D_2\}$ with   $D_1= \sup_t\{\delta_1(t)  \}$ and $D_2= \sup_t\{\delta_2(t)  \}$.\\
 Since $p>q>1$ we make use of the inequality
 \begin{equation*}
  \label{ineq}
\Psi_2^q = (a \Psi_2^p)^{\frac q p} \Big(a^{-\frac q{p-q} }  \Big)^{\frac {p-q }p}\leq 
\frac q p ( a \Psi_2^p) + \frac {p-q }p a^{-\frac q{p-q}  },
\end{equation*}
valid for abitrary $a>0$.\\
Choosing $a= \frac p q$ we arrive at
 \begin{equation}
  \label{ineq Psi_2^q}
 \Psi_2^p \geq \Psi_2^q - Q ,
\end{equation}
 with $Q:=  \frac {p-q }p  (\frac q p)^{\frac q{p-q}  }.$\\
 Inserting \eqref{ineq Psi_2^q} in \eqref{Psi'bis} we obtain the first order differential inequality
 \begin{equation}\label{Psi' ter}
\Psi' \geq  - \delta \Lambda_1\Psi  + c(\Psi_1^q + \Psi_2^q) - c\, Q \geq  -\delta \Lambda_1\Psi + 2^{1-q}c \Psi^q  - c\, Q := H(\Psi),
 \end{equation}
  where in the second term of \eqref{Psi' ter} we used the aritmetic inequality
\begin{equation}
\label{ineq 1}
 X^q + Y^q  \geq 2^{1-q} (X+Y)^q, \ q>1,
\end{equation}
with  $X \geq 0, \ Y \geq 0.$  
\\
  Since the initial data satisfy \eqref{H(Psi)}, then $\Psi(t)$ is increasing for small values of $t$. We have that $H(\Psi)$ is increasing in $\Psi$ from its negative minimum and it follows that $H(\Psi(t))$ is increasing for $t>0$. Thus $\Psi'(t)$ remains positive, so that $\Psi(t)$ blows up at time $t^*$.\\
 Integrating \eqref{Psi' ter} we obtain the upper bound $T_0$ in \eqref{upper}.
\qed}
\end{prth1.1}

In the case  $q=p$ we obtain an explicit estimate from above for $t^*$.  
\begin{proofcor}
{\rm 
We restart from the inequality \eqref{Psi'bis} and put $q=p$:
\begin{equation}\label{Psi' p=q}
\Psi' (t)  \geq   -\delta \Lambda_1 \Psi + c (\Psi_1^p + \Psi_2^p),
 \end{equation}
Since $p>1$, in  \eqref{Psi' p=q}, we apply the inequality  \eqref{ineq 1} with $q$ replaced by $p$.
Then
\begin{equation}\label{Psi' final p=q}
\Psi' (t)  \geq  -\delta \Lambda_1 \Psi + \bar c \Psi^p ,
 \end{equation}

with $\bar c= 2^{1-p} c.$
Integrating \eqref{Psi' final p=q} from $0$ to $t$, we arrive at
\begin{equation}
\label{integr Psi' final p=q}
\begin{aligned} 
\Psi^{1-p}(t)\leq  e^{(p-1) \delta \Lambda_1 t}   \Big( \Psi_0^{1-p} - \frac {  \bar c } { \delta \Lambda_1 }\Big) +  \frac {  \bar c } {\delta  \Lambda_1 }=: \mathcal{H}(t).
\end{aligned}
 \end{equation}
 Since the initial data satisfy \eqref{Psi0 p=q}
then $\mathcal{H} (t)$ vanishes at some time $\overline  T$. As a consequence, $\Psi(t)$ must blow up at some time $t^* \leq \overline  T$ with $\overline  T$ the desired upper bound \eqref{bar T, p=q}.}
\qed
\end{proofcor}
%%==========================================================%%
%%==============                                                    ==============%%
%%======                                 Section                                    ======%%
%%====                                                                                        ====%%
%%==                                                                                              ==%%
%%===                                    Lower Bound.                                ====%%
%%======                                                                                 ======%%
%%==============                                                     ==============%%
%%==========================================================%%

\section*{Acknowledgments}
M. Marras and S. Vernier-Piro are members of the Gruppo Nazionale per l'Analisi Matematica, la Probabilit$\grave{\rm a}$ e le loro Applicazioni (GNAMPA) of the Istituto Nazionale di Alta Matematica (INdAM). 

\subsection*{Financial disclosure}

M. Marras is partially supported by the research project: {\it Evolutive and stationary Partial Differential Equations with a focus on biomathematics} (Fondazione di Sardegna 2019), and by the grant PRIN n. PRIN-2017AYM8XW: {\it Non-linear Differential Problems via Variational, Topological and Set-valued Methods}.

Email address: mmarras@unica.
\\ 
Email address: svernier@unica.it

\end{document}